\documentclass[papersaver]{article}

\def\Spec{{ \mbox{Spec} }}
\def\crys{{ \mbox{crys} }}
\def\Hom{{ \mbox{Hom} }}
\def\MF{{ {\bf MF} }}
\def\1ox{{ \Omega^1_{\scriptstyle{X}} }}
\def\2ox{{ \Omega^2_{\scriptstyle{X}} }}

\def\ok1{{ \Omega^1_K }}
\def\ok2{{ \Omega^2_K }}

\def\ra{{ \rightarrow }}

\def\hra{{ \hookrightarrow }}

\def\A{{ {\bf A } }}

\def\8{{ {\infty } }}

\def\G{{ \Gamma }}
\def\Gal{{ \mbox{Gal} }}

\def\^{{ ^{\wedge} }}

\def\Z{{ {\bf Z } }}

\def\k{{ {\bar k} }}

\def\K{{ \bar{K} }}
\def\tor{{ \mbox{tor} }}

\def\MF{{ {\bf MF} }}
\def\Crys{{ \mbox{ Crys} }}
\newtheorem{thm}{Theorem}

\usepackage{amsmath}
\usepackage{amsfonts}
\usepackage{amscd}
\usepackage{amssymb}
\def\Kbar{\overline{K}}
\def\Rbar{\overline{R}}
\def\Qp{\mathbb{Q}_p}
\def\Zp{\mathbb{Z}_p}
\def\invlim{\varprojlim}
\def\dirlim{\varinjlim}

\title{A note on the crystalline subrepresentation functor}
\author{Minhyong Kim and Susan H. Marshall}

\begin{document}

\nocite{bures}

\maketitle\medskip

{\bf Abstract}
\smallskip

We propose the notion of the
{\em crystalline sub-representation functor} defined on
$p$-adic representations of the Galois groups
of finite extensions of $\Qp$, with certain restrictions
in the case of integral representations. By studying its
right-derived functors, we find a natural extension of
a  formula of Grothendieck expressing
the group of connected components of a Neron model
of an abelian variety in terms of Galois cohomology.

\medskip

\begin{center}
{\bf Une remarque sur le foncteur de sous-representation cristalline}
\end{center}

\medskip

{\bf Resum\'{e}}:
\smallskip

Nous proposons la notion d'un {\em foncteur de sous-representation
cristalline} defini pour les representations p-adiques
des groupes de Galois des extensions finies de $\Qp$,
avec certaines restrictions dans le cas des representations
integrales. Nous \'{e}tudions leur foncteurs deriv\'{e}s \`{a}
droite et les utilisons pour obtenir une
g\'{e}n\'{e}ralisation naturelle d'une formule de Grothendieck
donnant le groupe de composantes d'un mod\`{e}le de Neron
d'une vari\'{e}t\'{e} ab\'{e}lienne en terme de cohomologie galoisienne.

\medskip

Notation:

$K$: finite extension of $\Qp$.

$R$: the ring of integers in $K$.

$K_0$: maximal unramified subfield of $K$.

$k$: residue field of $K$=residue field of $K_0$.

$W=W(k)$: Witt vectors of $k$= ring of integers in $K_0$.

$\bar{K}$: an algebraic closure of  $K$.

$K^u$: the maximal
unramified subextension of $\K/K$.

$G=\Gal(\K/K)$: the Galois group of $\K$ over $K$.

$I=\Gal (\K/K^u)$: the inertia subgroup of $G$.

$l$: a prime different from $p$.

\vspace{5mm}

\section{ Maximal Crystalline Subrepresentations}
It is well known that the representations of $G$ with coefficients in
$\Zp$-modules, the {\em $p$-adic representations}, have very
different properties from the representations in $\Z_l$-modules for
$l\neq p$.  For example, even for a variety over $K$ with good
reduction over $R$, the representation of $G$ on the $p$-adic
\'{e}tale cohomology is only rarely unramified.

On the other hand, $p$-adic Hodge theory has provided us with a fine
classification of $p$-adic presentations together with appropriate
analogies to the $l$-adic case. For example, the $p$-adic notion
corresponding to an unramified $l$-adic representation is that of a
{\em crystalline} representation. These are the representations that
correspond via $p$-adic Hodge theory to weakly-admissible
 crystals (the correct $p$-adic
analogue of local systems), whereas representations that are genuinely
unramified correspond to the much smaller subcategory consisting of
crystals of slope zero (see, for example, \cite{F-icm}).

We wish to continue this analogy by presenting a new class of
cohomology theories associated to $p$-adic representations of Galois
groups of local fields. The definition is very natural and elementary,
and is likely to be well-known to experts. However, a specific
application motivated us to commit at least a short exposition to
paper:

Let $A$ be an abelian variety over $K$ and let $\A$ be its Neron model
over $R$. Let $A_0$ be the special fiber of $\A$ and $A_0^0$ the
connected component of the identity in $A_0$. Finally let
$\G=A_0(\k)/A_0^0(\k)$ be the geometric points of the group of
connected components of $A_0$.  Grothendieck points out the following
formula expressing the $l$-primary part of $\G$ in terms of Galois
cohomology:
$$\G(l)=H^1(I,T_l(A))_{\tor}$$ 
where $T_l$ refers to the $l$-adic Tate
module and the subscript denotes the torsion subgroup.  The motivating
problem is that of expressing the $p$ part of $\G$ in an analogous
`cohomological' manner involving only the generic fiber.

The formula is definitely false in general if we simply substitute $p$
for $l$. An easy argument using Kummer theory shows that when $A$ is
semi-stable over an absolutely unramified base, we actually have an
injection
$$H^1(I,T_p(A))_{\tor}\hra \G(p)$$ which is non-surjective in
general. For example, we can consider the case of an elliptic curve
with split semi-stable reduction and order of discriminant $p$.  It is
an easy exercise to check that in that case, the map is surjective iff
the elliptic curve has an unramified point of order $p$ which occurs
exactly when its Tate parameter is a $p$-power in $K^u$.  In short,
the torsion in the Galois cohomology of $I$ is not big enough to
capture the $p$-part of the component group.  But notice that the
Galois cohomology $H^1(I,\cdot)$ is just the first (right-)derived
functor of the functor
$$(\cdot) \mapsto (\cdot)^I$$ which we view as assigning to a
representation its {\em maximal unramified subrepresentation}. This is
an example of a `subrepresentation functor' or a `subobject' functor,
which can occur in a wide variety of contexts whenever one has
suitable subcategories of categories. On the other hand, we have
already remarked that the unramified objects comprise a sub-category
too small for geometric applications related to $p$-adic
representations.  This motivates us to define the {\em crystalline
subrepresentation functor} $$\Crys$$ from the category of
$\Qp$-representations of $G$ to itself.  Given a $\Qp$
representation $V$ of $G$, $\Crys(V)$ is the maximal crystalline
subrepresentation of $V$, where crystalline is defined in the usual
way for finite-dimensional representations and in general, we say $V$
is crystalline if it is a direct limit of finite-dimensional
crystalline subrepresentations.  Equivalently, we could say $V$ is
crystalline iff any finite dimensional subrepresentation is
crystalline. This equivalence follows from the fact that the category
of finite-dimensional crystalline representations is closed under
sub-objects.  The fact that it's also closed under quotient objects
implies that there is a well-defined notion of a `maximal' crystalline
subrepresentation.  The functor $\Crys$ is the natural $p$-adic
analogue of the `invariants under inertia' functor on $l$-adic
representations from the point of view of sub-representation
functors. Consequently, the derived functors of $\Crys$ are natural
analogues of Galois cohomology with respect to $I$.

To see that these notions are well-defined,
we must check two things:

(1) $\Crys$ is indeed a functor: This follows from
the fact that a quotient of a crystalline representation
is also crystalline, so that under a map
$V\ra W$ of representations, the crystalline part
must land in the crystalline part.

(2) $\Crys$ is left exact: The key point is that if
$U\subset V$ is a subrepresentation, then
$\Crys(U)=U\cap \Crys(V)$. The inclusion in the two
directions follows from the maximality involved in the
definition and the  sub-object property mentioned
earlier.

One could equally easily define the various `truncated' functors
$\Crys_{[a,b]}$ which associates to a representation the maximal
subrepresentation with Hodge-Tate weights in the interval $[a,b]$. We
will concentrate mostly on the functors $\Crys_{[0,h]}$ which we will
abbreviate as $\Crys_h$. It will be convenient to use the term {\em
$h$-crystalline representations} for the objects in the image of this
functor.  It is interesting to note that $\Crys_0$ is nothing but the
old inertia-invariants functor, so that the sequence of functors
$\Crys_0, \Crys_1, \ldots$ and their derived functors provide natural
prolongations of Galois cohomology. We see also that $\Crys$ is a bit
more than just an `analogue' of the inertia invariants
functor. Rather, the existence of these prolongations reflect the
richer structure that $p$-adic representations tend to have compared
to their $l$-adic counterparts. We propose that these derived functors
are natural invariants of $p$-adic representations (at least as
natural as Galois cohomology) and should be studied seriously. One
reason for thinking so stems from the application mentioned above.
For this, we need to define these functors also for integral
representations. Unfortunately, here the existing techniques for
making the correct definitions are rather incomplete, and we can
define {\em only} the truncated functors $\Crys_i$ for $i\leq
p-2$. (One can actually prolong it slightly to $i=(p-1)^*$ in an
appropriate sense, but we shall keep to the smaller truncation for
simplicity of exposition.)

We also need to assume that $K$ is absolutely unramified so that
$K=K_0$ and $R=W$. The foundational material we need is contained in
the seminal paper of Fontaine and Laffaille \cite{F-L}, but the reader
can find a nice summary in \cite{W}.

Let $h$ be a natural number $\leq p-2$.  One first defines finite
crystalline representations of height $\leq h$, or the finite {\em
$h$-crystalline representations }, to be the essential image of the
category $\MF^h_{R, \tor}$ (the finite-length filtered
$\phi$-modules of height $\leq h$) under the fully-faithful functor
$$M\mapsto V^*_{\crys}(M):=\Hom_{\MF_R}(M,A_{\crys,\infty})$$

Next, one defines a finite-type $\Zp$-module $L$ with $G$-action to
be $h-$crystalline if $L=\invlim L_i$ where the $L_i$ are
finite-length $h$-crystalline representations. The fact that
$h-$crystalline representations are closed under sub- and quotient
objects follows from the corresponding property for
$\MF_{R,\tor}^h$. In particular, this implies that a finite-type
$\Zp$ representation $L$ is $h$-crystalline iff $L/p^nL$ is
$h-$crystalline for all $n$ (which is the definition of \cite{N}), and
when $L$ is free, iff
$$L=V^*_{\crys}(M):=\Hom_{\MF_R}(M,A_{\crys})$$ for an object $M$ of
$\MF^h_R$  (the finitely generated
free filtered $\phi$-modules of height $\leq h$) (\cite{W} 2.2.2).

Now for an arbitrary $\Zp[G]$-module $V$, we define it to be
crystalline if $V=\dirlim V_i$ where the $V_i$ are
subrepresentations of finite-type.

We need to check that this definition is consistent with the existing
one for $\Qp$-representations.  Since we defined it for the
infinite-dimensional case using limits from finite dimensions, we need
only check it for finite-dimensional representations. So assume that
$V$ is $h-$crystalline in the old sense. Then
$V$=$\Hom_{MF_K}(\Delta,B_{\crys})$ for some $\Delta$ in $\MF_K^h$
\cite{F-L} (remarque 8.5 and 8.13 (c)).  Since $\Delta $ is
$B_{\crys}$-admissible, in particular, weakly admissible, one can find
a strongly divisible lattice $M\subset \Delta$ which is an object of
$\MF^h_R$. So we get $V=L\otimes \Qp$ where $L=\Hom_{\MF_R}(M,
A_{\crys})$. Now, $L$ is $h-$crystalline and $V=\dirlim L[1/p^n]$
while $L[1/p^n]\simeq L$ (via multiplication by $p^n$) is
$h$-crystalline.  So $V$ is $h$-crystalline in the new sense.  In the
other direction, assume $V=\dirlim L_i$ for $h-$crystalline
submodules $L_i$ of finite-type. Then some $L=L_i$ is a lattice and
$V=L\otimes \Qp$. But $L=\Hom_{\MF_R}(M, A_{\crys})$ for some free
$R$-module $M$ in $\MF_R^h$ and $M$ is then a strongly divisible
lattice in $\Delta:=M\otimes K$ according to the terminology
of \cite{F-L} definition 7.7, and therefore, $M\otimes K$ is weakly
admissible. Thus, by the main theorem of \cite{F-L}, $M\otimes K$ is
$B_{\crys}-$admissible and $V=\Hom_{\MF_K}(M\otimes K, B_{\crys})$ is
crystalline.

Thereby, we can define $\Crys_h$, the maximal $h-$crystalline
subrepresentation functor for $h\leq p-2$ compatibly on all $\Zp[G]$
modules. 
An easy consequence of the definitions is
that if $L$ is a finitely generated free $\Zp$ representation,
then $\Crys_h(L)=\invlim \Crys_h(L/p^nL)$.
It should be emphasized that we also have $\Crys$ and all the
other $\Crys_{[a,b]}$'s if we stick to rational representations.  By
the key property that $$L_1 \subset L_2\Rightarrow
\Crys_h(L_1)=L_1\cap \Crys_h(L_2)$$ we again have left exactness and
therefore, all the right-derived functors.

A systematic study of these functors will be presented in the
forthcoming Ph.D. thesis of the second author.

\section{The $p$-complement to Grothendieck's formula}

In this section, we will continue to assume that $K$ is
absolutely unramified, and furthermore, that
$p>2$.

We will be using one more functor $FF$ which associates to a $p$-adic
representation its maximal `finite and flat' part.  Of course, one
needs to define finite flat $p$-adic representations in a general
setting. For finite $\Zp$ representations, finite flat means the
usual thing: a finite representation is finite flat if it's isomorphic
to the $\K$ points of a finite flat commutative
group scheme over $R$. A
finite-type $\Zp$-representation is defined to be finite flat if it
is the inverse limit of finite finite flat representations (the double
adjective seems unfortunately unavoidable). Finally, an arbitrary
$\Zp$ representation for $G$ is said to be finite flat if it is the
direct limit of finite flat representations of finite type.

For finite representations, the property of being finite flat is
closed under passing to sub-objects and quotient objects (using
Zariski closure and construction of good quotient schemes), so the same
is true for any $\Zp$ representation. Thus it makes sense to speak of
the maximal finite flat subrepresentation of any representation, and
the associated functor $FF$ is left exact. Thus, we can consider its
derived functors.  In fact, by 
Fontaine-Laffaille's description of finite
flat group schemes (\cite{F-L}, section 9) $FF$ is nothing but $\Crys_1$.
Notice, however, that $FF$ is defined over an arbitrary
local field, not necessarily absolutely unramified.

We will also need the trivial observation that if
$$0 \ra M_1 \ra M_2 \ra M_3 \ra 0$$ is an exact sequence in $\MF_R$,
$M_1$ and $M_3$ are in $\MF_R^h$, and $M_2$ is in $\MF_R^{h'}$ for
some $h'$, then in fact, $M_2$ is in $\MF_R^h$.  This follows by
noting that the morphisms are strict so that any $F^iM_2$ for $i>h$
would have to be zero when intersected with $M_1$ and mapped to $M_3$,
and hence, must be zero.

Thus, we have an obvious corresponding statement for $h-$ and
$h'$-crystalline representations.

We now return to the problem of expressing the $p$-part of $\Gamma$ in
an analogous manner to Grothendieck's formula for $l \neq p$
\begin{equation}
\Gamma(l) \cong H^1(I, T_lA)_{\tor}. 
\end{equation}
To derive the above formula, Grothendieck shows that 
\begin{equation}
\Gamma[l^n] \cong A[l^n]^f/(A^0[l^n])^f
\end{equation}
where $\Gamma[l^n]$ (resp. $A[l^n]$) denotes the
kernel of multiplication by $l^n$ on $\Gamma$ (resp. $A(\Kbar)$), 
 and the
superscript $f$ denotes the ``finite part'' (denoted the ``fixed
part'' by Grothendieck in \cite{sga7}, section 2.2.3), i.e.,
the points that extend to a map from $\Spec( \Rbar)$ to
$\A$, or equivalently, the $\Kbar$ points of the
maximal finite flat subgroup scheme of $\A[l^n]$.  
Similarly, $(A^0[l^n])^f$ denotes the 
 $\Kbar$ points of the maximal finite flat
subgroup scheme of $\A^0[l^n]$ which can also be thought
of as the points of $A[l^n]^f$
which reduce mod $p$ to a point in $A_0^0$, the connected
component of the identity in the special fiber.
The key point
then is that the finite part coincides with the inertia invariants of
$A[l^n]$ (resp. $A^0[l^n]$) \cite{sga7} (Proposition 2.2.5) and the
formula (1) follows easily.

In the case of $l=p$, (2) still holds (provided one assumes
semi-stability), but it is no longer the case that the fixed part and
inertia invariants coincide. However, we will show below that (for $l
\neq p$) the finite part coincides with the \emph{maximal $h$-crystalline
part} for any $1\leq h\leq p-2$ (recall that $p>2$).  
This will allow us to derive, in a completely analogous manner
to Grothendieck, the following:
\begin{thm} Let $A$ be an abelian variety over the
absolutely unramified local field $K$ with semi-stable
reduction and $1\leq h\leq p-2$. Then
$$\Gamma(p) \cong R^1\Crys_h(T_pA)_{\tor}.$$
\end{thm}

{\em Proof.}

We will first show that $\Crys_h(A[p^n]) = (A[p^n])^f$.
For this, we  note
that the fixed part of $A[p^n]$ is none other than $FF(A[p^n])$. That is,
 the
fixed part is finite-flat by definition, giving us one inclusion 
\[ (A[p^n])^f \subset FF(A[p^n]). \]
Now let $\mathcal{V}$ denote the finite-flat group scheme extending
$FF(A[p^n])$, so that if $V$ is the generic fiber of $\mathcal{V}$, we have
\[ V(\Kbar) \cong FF(A[p^n]) \] as $G$-modules. From the inclusion 
$FF(A[p^n]) \subset A[p^n]$, we have a map
\[ V \rightarrow A. \] We need to show that this map
extends to a map $\mathcal{V} \ra \A$, thereby showing that
the finite part is actually ``finite inside $\A$.''
However,
restricting to the connected component $V^0$ of $V$, we find that the
image must actually land in the finite part of $A$. This follows
because $A[p^n]/A[p^n]^f$ is unramified (\cite{sga7}, Proposition 5.6).
  By
results of Raynaud \cite{R}, this extends to a map $\mathcal{V}^0 \rightarrow
\A^f$. Hence by Lemma 5.9.2 of \cite{sga7}, we get a unique map
$\mathcal{V} \rightarrow \A$ extending the two previous maps,
and giving us the opposite inclusion.
(This is essentially the same argument as in \cite{ribet}, Lemma 6.2.)

We saw above that $FF = \Crys_1$, as functors.  We will now
show that one can replace $\Crys_1$ by any
of the $\Crys_h$'s in our setting. In fact,
we will see from the proof that if any general $\Crys$ functor
were defined for finite representations, then that
could be used as well.

We certainly have an
inclusion
\[ FF(A[p^n]) \hookrightarrow \Crys_h(A[p^n]) \]
which induces an inclusion of $\Crys_h(A[p^n])/FF(A[p^n])$ into the
unramified $G$-module $A[p^n]/FF(A[p^n])$. Thus,
$\Crys_h(A[p^n])/FF(A[p^n])$ is unramified as well, and hence
finite-flat as a representation. Therefore $\Crys_h(A[p^n])$ sits in
the middle of a short exact sequence whose outer terms are both
crystalline of height one (actually, the last is of height 0):
\begin{equation*}
0 \rightarrow FF(A[p^n]) \rightarrow  \Crys_h(A[p^n])  \rightarrow 
	 \Crys_h(A[p^n])/FF(A[p^n])   \rightarrow 0 
\end{equation*}
By the  observation made earlier, we see that $\Crys_h(A[p^n])$ is itself
crystalline of height one, and thus equal to $\Crys_1(A[p^n])=FF(A[p^n])
=(A[p^n])^f$. As $A^0[p^n]/FF(A^0[p^n])$ is contained in
$A[p^n]/FF(A[p^n])$, it is also unramified and an entirely similar argument
gives $\Crys_h(A^0[p^n]) = (A^0[p^n])^f$ as well.

The isomorphism (2) thus becomes
\[ \Gamma[p^n] \cong \Crys_h(A[p^n])/\Crys_h(A^0[p^n]) \]
However,
\vspace{3mm}

{\bf Claim:}
\[\Crys_h(A[p^n])/\Crys_h(A^0[p^n]) \simeq 
\Crys_h(T_pA \otimes \Zp/p^n\Zp)/\Crys_h(T_pA) \otimes \Zp/p^n\Zp. \]

{\em Proof.} The equality between the `numerators' is
obvious, so we need to
 see that $\Crys_h(T_pA) \otimes \Zp/p^n\Zp$ is equal to
$\Crys_h(A^0[p^n])$. But $$\Crys_h(T_pA) \cong
\invlim(\Crys_h(A[p^n])=\invlim(A[p^n]^f)$$ 
Hence  $\Crys_h(T_pA)\simeq (T_pA)^f$. Thus, 
\[ \Crys_h(T_pA) \otimes \Zp/p^n\Zp \cong (T_pA)^f \otimes \Zp/p^n\Zp
\cong (A^0[p^n])^f \cong \Crys_h(A^0[p^n]). \] 
(The key point is the second isomorphism, as explained in
\cite{sga7}. That is, if you take the finite part of 
the Tate module and
then reduce mod $p^n$, then you end up in $(A^0[p^n])^f$
because the multiplication by $p$ map is finite
and surjective only on $\A^0$.)
\vspace{3mm}

Applying direct limits, we get the formula
\[ \Gamma(p) \cong \Crys_h(T_pA \otimes \Qp/\Zp)/\Crys_h(T_pA) \otimes \
\Qp/\Zp. \] This plays a role analogous to Grothendieck's formula
\cite{sga7} (Proposition 11.2).

Following Grothendieck, we next apply $\Crys_h$ to the short exact sequence
\begin{equation*}
\begin{CD}
0 @>>> T_pA  @>>> T_pA \otimes \Qp @>>> T_pA \otimes \Qp/\Zp  @>>> 0 
\end{CD}
\end{equation*}
and obtain the long exact sequence
\begin{equation*}
\begin{CD}
0 @>>> \Crys_h(T_pA)  @>>> \Crys_h(T_pA \otimes \Qp) @>>> 
	\Crys_h(T_pA \otimes \Qp/\Zp)  @>>> \\
  @>>>  R^1\Crys_h(T_pA)  @>>>  R^1\Crys_h(T_pA \otimes \Qp) @>>>  
	R^1\Crys_h(T_pA \otimes \Qp/\Zp)  @>>> \dots 
\end{CD}
\end{equation*}
Note that $$\Crys_h(T_pA \otimes \Qp) = 
(\Crys_h(T_pA \otimes \Qp)\cap T_pA)\otimes \Qp=\Crys_h(T_pA) \otimes \Qp$$
  Thus the kernel of the map of 
\[ R^1\Crys_h(T_pA) \rightarrow R^1\Crys_h(T_pA \otimes \Qp) \]
 is $\Crys_h(T_pA \otimes \Qp/\Zp)/\Crys_h(T_pA) \otimes \Qp/\Zp$, i.e.
$\Gamma(p)$. Since $\Gamma(p)$ is torsion and $R^1\Crys_h(T_pA \otimes
\Qp)$ torsion-free, we do indeed find that
\[ \Gamma(p) \cong R^1\Crys_h(T_pA)_{\tor}. \]

\vspace{3mm}

{\bf Remark.}
From the proof, it is clear that one could have just
used the functor $FF$ for 
the theorem in which case one could extend the theorem
to the case of $e\leq p-2$ by eliminating the
Fontaine-Laffaille theory. However, this would have
made the analogy to the $l$-case less natural,
since a crystalline resepresentation
is clearly the correct general notion which 
sets the formula into a broad context. In particular,
the definition of $FF$ on $\Qp$ representations
is rather artificial compared to $\Crys$.
It would of course have been nicer to replace
$\Crys_h$ by a general $\Crys$ even for the
integral representations.

\vspace{3mm}
{\bf Acknowledgements:} We are grateful to Wiesia Niziol for
pointing out the definitions in \cite{N} and to Ken Ribet
for directing us to lemma 6.2 of \cite{ribet}.
We are especially grateful to
Kirti Joshi for innumerable conversations on p-adic Hodge
theory. 

M.K. was supported in part by NSF grant DMS-9701489

{\footnotesize DEPARTMENT OF MATHEMATICS, UNIVERSITY OF ARIZONA,TUCSON, AZ 85721, U.S.A. E-MAIL: kim@math.arizona.edu susan@math.arizona.edu}
\end{document}